\documentclass[12pt]{article}
\usepackage{amsfonts,amstext}
\begin{document}

\newtheorem{Th}{Theorem}[section]
\newtheorem{Def}{Definition}[section]
\newtheorem{Lemm}{Lemma}[section]
\newtheorem{Cor}{Corollary}[section]
\newtheorem{Rem}{Remark}[section]

\title{\bf Chern numbers of Chern submanifolds}
\author{K.~E.~Feldman\footnote{e-mail:feldman@maths.ed.ac.uk}}
\date{}
\maketitle

\begin{center}
Department of Mathematics and Statistics, University
of Edinburgh,\\ James Clerk Maxwell Building,  Mayfield 
Road, Edinburgh, Scotland,\\ EH9 3JZ 
\end{center}

\begin{abstract}
We present a solution of the generalized Hirzebruch problem on the 
relations between the Chern numbers of a stably almost complex 
manifold and the Chern numbers of its virtual Chern submanifolds.
\end{abstract}

\section{Introduction}
 A cobordism class $\alpha\in \Omega^k(M)$ of a smooth 
manifold $M$ can be realized by
a smooth submanifold $L\subset M\times {\mathbb R}^N$ of codimension $k$
for $N$ large enough and $L$ is unique up to bordism.  We  call
such a submanifold a virtual submanifold of $M$.
 For example, the virtual submanifold corresponding
to the Euler class $\chi(\xi)$ of a vector bundle $\xi$ is 
realized by the zero set of any generic section $s:M\to \xi$. 
If the cobordism theory $\Omega^*(\cdot)$ 
has a special structure then the normal 
bundle of the virtual submanifold $L$ also has the corresponding structure.
It is natural 
to expect that virtual submanifolds inherit some properties from the
original manifold.
 In~\cite{BV} the authors found a non-trivial relation between the signatures 
of an almost complex manifold $M^{4n}$ and the virtual submanifolds 
$[P_k(\tau(M^{4n}))]$ dual to the complex cobordism  Pontrjagin 
classes $P_k(\tau)\in U^{4k}(M^{4n})$~\cite{CF}.
They posed a more general question: What are all the divisibility
relations between the Chern numbers of a stably almost complex manifold and
its virtual tangent Chern submanifolds? This is a generalization of the 
classical Hirzebruch problem~\cite{Hi1} on the divisibility restrictions for 
the Chern numbers of a stably almost complex manifold, which was solved by 
Hattori and Stong~\cite{Ht,St}.

In the present paper we find all the relations between the Chern 
numbers of a stably 
almost complex manifold $M^{2n}$ and its virtual submanifolds $[c_k(\eta)]$
defined by the
Chern classes $c_k(\eta)\in U^{2k}(M^{2n})$ of 
a complex vector bundle $\eta$. The result 
of~\cite{BV} is a particular case of our general theorem.  

For the solution we construct a multiplicative transformation $\tau$ of 
complex cobordism  theory $U^*(\cdot)$ to ordinary cohomology 
theory $H^*(\cdot,\Lambda[y]\otimes\Lambda[z]\otimes \mathbb Q)$,
where following~\cite{Mac} $\Lambda[y],\Lambda[z]$ are the rings of symmetric 
polynomials over $\mathbb Z$
in variables $y_i$ and $z_i$, $i=1,2,...$, respectively. 
The key tools for the construction are the 
Chern--Dold character~\cite{B} and the Landweber--Novikov 
operations~\cite{L,N}. 
Applying the 
Riemann--Roch theorem~\cite{AH1} to the transformation $\tau$  we can express 
any Chern number of the virtual Chern 
submanifold in terms of the Chern numbers of the original manifold.  
Together 
with the Hattori--Stong 
theorem this expression  describes all the relations 
between Chern numbers.

 The paper is organized in the following way. 
In Section 2 we define the transformation 
$\tau:U^*(\cdot)\to H^*(\cdot;\Lambda[y]\otimes \Lambda[z]\otimes {\mathbb Q})$
and calculate its Todd genus.

 In Section 3 we apply the Riemann--Roch theorem to the transformation $\tau$
and obtain an expression for the normal Chern numbers of  
the normal virtual Chern submanifolds $[c_k(\nu(M^{2n}))]$ in terms of 
the Chern numbers of the original manifold $M^{2n}$.

 In Section 4 we present some examples of numerical relations 
between Chern numbers of Chern submanifolds. In particular, we deduce the 
formula obtained by Buchstaber and Veselov in~\cite{BV}.

 In Section 5 we prove a general theorem on the relations between Chern 
numbers of 
virtual Chern submanifolds $[c_k(\eta)]$ corresponding to a complex vector 
bundle $\eta$ over a stably almost complex manifold $M$. Using the 
Hattori--Stong theorem we deduce a new divisibility condition for 
the cohomology Chern classes. As an example we prove that the 
Euler 
characteristic of an almost  complex manifold whose tangent bundle possesses 
a complex line subbundle is even.
 
\section{The transformation $\tau$ and its Todd class}

Denote the Landweber--Novikov operation corresponding to the symmetric 
polynomial
$f\in \Lambda[z]$ by ${\bf f}$. For the Thom class $u_m\in U^{2m}(MU(m))$ 
and the tautological bundle $\eta_m$ over $BU(m)$ we have
$$
{\bf f}(u_m)=u_m f(\eta_m),
$$
where $f(\eta_m)$ is the Chern class in complex cobordism 
corresponding to 
the symmetric polynomial 
$f$. Let $p_k(z)\in \Lambda [z]$ be the $k$-th power sum of $z_i$, $i=1,2,...$ 
Denote the monomial symmetric function corresponding to the partition 
$\lambda=(\lambda_1\ge\lambda_2\ge...\ge 0)$ by 
$m_{\lambda}(z)$. One defines the complete symmetric functions $h_{\lambda}(z)$
using the identity~\cite{Mac}
\begin{equation}
\label{SymId}
\prod^{\infty}_{i,j=1}(1-z_iy_j)^{-1}=
\exp\left(\sum^{\infty}_{k=1}\frac{p_k(z)p_k(y)}{k}\right)=
1+\sum_{\{\lambda||\lambda|>0\}}h_{\lambda}(z)m_{\lambda}(y).
\end{equation}
We combine all the Landweber--Novikov operations into one multiplicative 
operation ${\cal S}_{(z)}$:
$$
{\cal S}_{(z)}=\exp\left(\sum^{\infty}_{k=1}\frac{p_k(z){\bf p_k}}{k}\right)=
\sum_{\lambda} h_{\lambda}(z){\bf m_{\lambda}}.
$$ 
The operation ${\cal S}_{(z)}$ acts as a ring 
endomorphism on $U^*(X)\otimes \Lambda[z]$ for any topological space $X$. Moreover, this 
operation is invertible (see~\cite{BBNU} for further algebraic
properties of ${\cal S}_{(z)}$). 

In~\cite{B} the Chern--Dold character $ch_U:U^*(X)\to H^*(X,\Omega^U_*\otimes {\mathbb Q})$ was constructed 
which is the canonical embedding for $X=\{pt\}$. The value of $ch_U$ on the Thom class 
$u\in U^2(MU(1))$ is given by the formula~\cite{B}:
$$
ch_U(u)=x+\sum^{\infty}_{n=1}[N^{2n}]\frac{x^{n+1}}{(n+1)!},
$$
where $x$ is the Thom class of $MU(1)$ in ordinary cohomology, $N^{2n}$ is a 
closed stably almost complex manifold defined uniquely by 
the conditions $m_{\lambda}(N^{2n})=0$ for $\lambda\neq (n)$ and  
$m_{(n)}(N^{2n})=(n+1)!$
(we denote the normal Chern number of a stably almost complex manifold
$M^{2n}$ corresponding to the normal Chern class 
$f(\nu(M^{2n}))$, $\deg f=n$, by 
$f(M^{2n})$).

Define the ``universal'' symmetric genus
$$
S^*_{(y)}:\Omega^U_*\to  \Lambda[y]\qquad {\mathrm {by}}\qquad 
S^*_{(y)}([M^{2n}])=\sum_{\{\lambda||\lambda|=n\}}h_{\lambda}(y)
m_{\lambda}(M^{2n}).
$$
From~(\ref{SymId}) it follows that the generating function of the universal symmetric genus $S^*_{(y)}$ is
$$
S^*_{(y)}(x)=\prod^{\infty}_{i=1}(1-y_ix)^{-1}.
$$
The transformation
$$
H^*(X,\Omega^U_*\otimes \Lambda[z]\otimes {\mathbb Q})\to 
H^*(X, \Lambda[y]\otimes \Lambda[z]\otimes {\mathbb Q})
$$
determined through the coefficient ring homomorphism 
$S^*_{(y)}\otimes id\otimes id:\Omega^U_*\otimes \Lambda[z]\otimes 
{\mathbb Q}\to 
\Lambda[y]\otimes \Lambda[z]\otimes {\mathbb Q}$
is also denoted by $S^*_{(y)}$.

Let $\tau$ be the composition $S^*_{(y)}\circ ch_U\circ {\cal S}_{(z)}$.
Obviously, $\tau$ is a multiplicative transformation of cohomology 
theories
$$
\tau:U^*(X)\to H^*(X,\Lambda[y]\otimes \Lambda[z]
\otimes {\mathbb Q}).
$$
Recall that for any multiplicative transformation $\mu:h^*_1(\cdot)\to h^*_2(\cdot)$ 
of complex oriented cohomology theories $h^*_1(\cdot),h^*_2(\cdot)$ the Todd class
$T_\mu$ is defined via the Riemann--Roch theorem~\cite{AH1}:
$$
p^{h_2}_!\left(\mu(x)\cdot T_{\mu}(\nu)\right)=\mu\left(p^{h_1}_!(x)\right),
$$
where $p:M_1\to M_2$ is a map of smooth manifolds with the stably complex oriented normal bundle 
$\nu\cong p^*(\tau(M_2))-\tau(M_1)$, $p^{h_1}_!$ and $p^{h_2}_!$ are the corresponding Gysin maps in
$h^*_1(\cdot)$ and $h^*_2(\cdot)$ respectively, $x\in h^*_1(M_1)$.
\begin{Th}
\label{Todd}
The Todd class corresponding to $\tau$ is
$$
T_{\tau}(\eta)=\prod^n_{i=1}\left(\prod^{\infty}_{j=1}(1-y_jx_i)^{-1}\cdot
\prod^{\infty}_{k=1}\left(1-\frac{z_kx_i}{\prod^{\infty}_{j=1}
(1-y_jx_i)}\right)^{-1}\right),
$$
here $n=\dim_{\mathbb C}\eta$, $x_i$, $i=1,...,n$, are the 
Chern roots in integral
 cohomology of the complex vector bundle $\eta$. 
\end{Th}
{\bf Proof.} By splitting principle it is sufficient to verify this formula for the case $n=1$.
For the Thom class $u\in U^2(MU(1))$ we obtain
\begin{equation}
\label{T1}
{\cal S}_{(z)}(u)=u\cdot 
\exp\left(\sum^{\infty}_{k=1}\frac{p_k(z)u^k}{k}\right)=
u\prod^{\infty}_{k=1}(1-z_ku)^{-1};
\end{equation}
For the manifolds $N^{2n}$ in the definition of the Chern-Dold character one has
$S^*_{(y)}(N^{2n})=(n+1)!h_{(n)}(y)$. Then (see also~\cite{Mac})
\begin{equation}
\label{T2}
S^*_{(y)}\circ ch_U (u)=x+\sum^{\infty}_{n=1}h_{(n)}(y)x^{n+1}=x\prod^{\infty}_{j=1}(1-y_jx)^{-1}.
\end{equation}
The statement in the theorem now follows using the 
multiplicative property of $S^*_{(y)}\circ ch_U$.

\begin{Rem}
\label{reason}
The multiplicative operation ${\cal S}_{(z)}$ transforms a stably almost
complex manifold into the collection of all Chern submanifolds.
The Chern--Dold character maps them into the coefficient ring $\Omega^U_*$.
The operation $S^*_{(y)}$ calculates the Chern numbers of the elements 
of $\Omega^U_*$.
\end{Rem}
 
\section{Riemann--Roch theorem for $\tau$}

Consider a smooth compact closed connected manifold $M^{2n}$ without boundary 
with fixed almost complex structure on its stable normal bundle 
$\nu(M^{2n})$. Every Chern class of the normal bundle 
$\nu(M^{2n})$ can be represented by a submanifold of $M^{2n}\times 
{\mathbb R}^{N}$ with 
almost complex normal structure  for $N$ large enough. 
We will denote the manifold corresponding 
to the characteristic class $f(\nu(M^{2n}))\in U^*(M^{2n})$ by $[f(\nu(M^{2n}))]$.
Let $\langle x,[M^{2n}]\rangle$ be the 
Kronecker product of 
$x\in H^{2n}(M^{2n},\Lambda[y]\otimes\Lambda[z]\otimes {\mathbb Q})$ with the
fundamental class 
$[M^{2n}]\in H_{2n}(M^{2n},\Lambda[y]\otimes\Lambda[z]\otimes {\mathbb Q})$;
if $\deg x\neq 2n$ we put $\langle x,[M^{2n}]\rangle=0$.

\begin{Th}
\label{relations} The following equality holds:
\begin{equation}
\label{rel}
\langle T_{\tau}(\nu(M^{2n})),[M^{2n}]\rangle=
S^*_{(y)}(M^{2n})+\sum_{\{\lambda||\lambda|>0\}}h_{\lambda}(z)
S^*_{(y)}([m_{\lambda}(\nu(M^{2n}))]).
\end{equation}
\end{Th}
\begin{Rem}
The left side of this equality depends only on the 
Chern numbers of the manifold
$M^{2n}$. The right side depends on the Chern numbers of 
the virtual Chern submanifolds. Moreover, because of the linear independence
of $h_{\lambda}(z)$ in $\Lambda[z]$ and $m_{\lambda}(y)$ in $\Lambda[y]$,
the Chern numbers of a Chern submanifold can be expressed in terms of the 
Chern numbers of $M^{2n}$  using~(\ref{rel}). 
\end{Rem}
{\bf Proof.} Denote the constant map $M^{2n}\to \{pt\}$ by $p$.
The corresponding Gysin map in complex cobordism is denoted by $p^U_!$
and that in cohomology by $p^H_!$. 
There is a canonical cobordism class 
$\sigma\in U^0(M^{2n})$ defined by 
the identity map $M^{2n}\to M^{2n}$.
It is obvious that 
\begin{equation}
\label{RHS}
p^U_!(\sigma)=[M^{2n},\nu(M^{2n})]\in \Omega^{-2n}_U(pt)=\Omega^U_{2n}.
\end{equation}
We apply the Riemann--Roch theorem in the form given in~\cite{AH1} to the 
transformation $\tau$ and the element $\sigma$. Using~(\ref{RHS}) we obtain
$$
p^H_!(\tau(\sigma)\cdot T_{\tau}(\nu(M^{2n})))=\tau([M^{2n},\nu(M^{2n})]).
$$
Because the normal bundle of the identity map $M^{2n}\to M^{2n}$ is trivial,
we have $\tau(\sigma)=S^*_{(y)}\circ ch_U(\sigma)$. 
According to the results of \cite{B},
$$
(S^*_{(y)}\circ ch_U)(\sigma)=1.
$$
Thus, 
$$
p^H_!(\tau(\sigma)\cdot T_{\tau}(\nu(M^{2n})))=
\langle T_{\tau}(\nu(M^{2n})),[M^{2n}]\rangle.
$$
On the other hand, the normal bundle of the constant map $M^{2n}\to \{pt\}$
is $\nu(M^{2n})$. Thus,
$$
{\cal S}_{(z)}([M^{2n},\nu(M^{2n})])=
[M^{2n}]+\sum_{\{\lambda||\lambda|>0\}}
h_{\lambda}(z)[m_{\lambda}(\nu(M^{2n}))]\in \Omega^U_*.
$$
Note that $ch_U$ is the identity on the coefficient ring $\Omega^U_*$. 
Applying $S^*_{(y)}$ we obtain Equation~(\ref{rel}).

\begin{Rem}
\label{cobord} From Theorem~\ref{relations} it follows that  
the Chern numbers of any virtual Chern submanifold of a
stably almost complex manifold depend only on the cobordism class of the
original manifold. This is to be expected because the composition of any two 
Landweber--Novikov operations is again a Landweber--Novikov operation.
\end{Rem}

\section{Examples}

Various substitutions into~(\ref{rel}) lead us to families of 
relations between Chern numbers of stably almost complex manifold
and its Chern submanifolds. Let us consider the Hirzebruch genus 
$\chi_y$~\cite{Hi} which arises from the power series
$$
Q(x)=\frac{x(1+y{\mathrm e}^{-x(1+y)})}{1-{\mathrm e}^{-x(1+y)}}.
$$ 

\begin{Th}
\label{Chi_y}For the tangent virtual Chern submanifolds $[c_k(\tau(M^{2n}))]$
the following relation holds
\begin{equation}
\label{Chi1}
(1+y)^n T(M^{2n})=\chi_y(M^{2n})+\sum^n_{k=1} y^k\chi_y([c_k(\tau(M^{2n}))]),
\end{equation}
where $T(M^{2n})$ is the classical Todd genus of the manifold $M^{2n}$.
\end{Th}
{\it Note:} the Chern numbers used in this theorem are those 
of the tangent bundle.

\

\noindent
{\bf Proof.}
Substitute $z_1=-y$, $z_2=z_3=...=0$ into~(\ref{rel}) and choose the values of
$y_k$ so that 
$$
S^*_{(y)}(x)=\prod^{\infty}_{j=1}(1-y_jx)^{-1}=Q(x)^{-1}.
$$
If we calculate the normal Chern number corresponding to the 
$\frac{1}{Q}$-genus of any 
stably almost complex manifold, we obtain its $\chi_y$ genus. Because 
$c_k(\tau(M^{2n}))=(-1)^k h_k(\nu(M^{2n}))$, for this choice of 
variables the right hand side of~(\ref{rel}) becomes
$$
\chi_y(M^{2n})+\sum^n_{k=1} y^k\chi_y([c_k(\tau(M^{2n}))]).
$$
Now we calculate the left hand side of~(\ref{rel}) for this choice of 
variables. We have the following equality
$$
Q(x)^{-1}\left(1+\frac{yx}{Q(x)}\right)^{-1}=
(Q(x)+yx)^{-1}=\frac{1-{\mathrm e}^{-x(1+y)}}{x(1+y)}.
$$
The corresponding normal Chern number of $M^{2n}$ is exactly 
the classical Todd genus of $M^{2n}$ multiplied by $(1+y)^n$. Thus, 
Equation~(\ref{Chi1}) follows.

\begin{Rem}
\label{Sub_y}
 Putting $y=-1$ or $1$ into~(\ref{Chi1}) we obtain equations on 
the Euler characteristics and the signatures of the tangent
virtual Chern submanifolds.  
\end{Rem}
In the same manner we can obtain relations between the Chern numbers of a stably 
almost complex manifold and its virtual Pontrjagin submanifolds.
Recall that the tangent virtual Pontrjagin submanifold $[P_k(\tau(M^{4n}))]$ is 
defined as $[(-1)^kc_{2k}({\mathbb C}\otimes\tau(M^{4k}))].$

\begin{Th}
\label{Sig}
For the signatures of the submanifolds 
$[P_k(\tau(M^{4n}))]$ of a stably almost complex manifold $M^{4n}$
the following relation holds
\begin{equation}
\label{sign}
2^{4n}\hat A(M^{4n})=\sigma(M^{4n})+\sum^n_{k=1} (-1)^k
\sigma([P_k(\tau(M^{4n}))]),
\end{equation}
where $\hat A(M)$ is the $\hat A$-genus of $M$.
\end{Th}
{\bf Proof.}
Substitute $z_1=-z_2=1$, 
$z_3=z_4=...=0$ into~(\ref{rel}) and choose the values of
$y_k$ so that 
$$
S^*_{(y)}(x)=\prod^{\infty}_{j=1}(1-y_jx)^{-1}=\frac{tanh(x)}{x}=\tilde L(x).
$$
If we calculate the 
normal Chern number corresponding to the $\tilde L$-genus of any 
stably almost complex manifold, we obtain its signature. 
For a stably almost complex manifold $M^{4n}$  
the Pontrjagin classes can be expressed in terms of the Chern roots $r_j$,
$j=1,...,2n$, of the tangent bundle:
$$
\sum^n_{k=0} P_k t^{2k} =  \prod^{2n}_{j=1} (1+ir_jt)(1-ir_jt),
$$
where $t$ is a formal parameter. 
Using again the fact
$c_k(\tau(M^{4n}))=(-1)^k h_k(\nu(M^{4n}))$, for this choice of 
variables, the right hand side of~(\ref{rel}) becomes
$$
\sigma(M^{4n})+\sum^n_{k=1} (-1)^k\sigma([P_k(\tau(M^{4n}))]).
$$
Now we calculate the left hand side of~(\ref{rel}) for this case. 
We obtain the following  equality
$$
\tilde L(x)\left(1-x^2\tilde L(x)^2\right)^{-1}=
\left(\tilde L(x)^{-1}-x^2\tilde L(x)\right)^{-1}=\frac{{\mathrm e}^{2x}-{\mathrm e}^{-2x}}
{4x}=\hat A(4x)^{-1}.
$$
The corresponding normal Chern number of $M^{4n}$ is exactly the
$\hat A$-genus of $M^{4n}$ multiplied by $2^{4n}$. Thus, 
Equation~(\ref{sign}) 
follows.

\begin{Rem}
Denote the number of ones in the binary expansion of the number $n$ by 
$\alpha(n)$. It is well known that the $\hat A$-genus of a stably almost 
complex 
manifold $M^{4n}$ multiplied by $2^{4n-\alpha(n)}$ is an integer~\cite{AH2}. 
Thus, from Theorem~\ref{Sig} we can deduce the main result of~\cite{BV}
$$
\sigma(M^{4n})+\sum^n_{k=1} (-1)^k\sigma([P_k(\tau(M^{4n}))])=0(\bmod 2^{\alpha(n)}).
$$
\end{Rem}

\section{Chern numbers of general Chern submanifolds}

Virtual Chern submanifolds $[m_{\lambda}(\eta)]$ can be defined for any stably almost complex bundle $\eta$ over 
a stably almost complex manifold $M$ with the fixed complex structure in the normal bundle $\nu=\nu(M)$. It turns 
out that the Chern numbers of such submanifolds can also be expressed in terms of the Chern numbers of $M$.  
Introduce two power series with coefficients in $\Lambda[y]\otimes\Lambda[z]$ and $\Lambda[y]$ respectively:
$$
T_1(x)=\prod^{\infty}_{k=1}\left(1-\frac{z_k x}{\prod^{\infty}_{j=1}(1-y_jx)}\right)^{-1}, \qquad{} T_2(x)=\prod^{\infty}_{i=1}(1-y_ix)^{-1}.
$$
For a complex vector bundle $\eta$, $\dim_{\mathbb C}\eta=n$, over a base space $X$ we define two non-homogeneous elements in $H^*(X,\Lambda[y]\otimes\Lambda[z])$ by the formulae
$$
T_1(\eta)=\prod^n_{i=1} T_1(x_i),\qquad{} T_2(\eta)=\prod^n_{i=1} T_2(x_i),
$$
where $x_i$ are the Chern roots of $\eta$ in ordinary cohomology. 

\begin{Th} 
\label{ChCh}
The Chern numbers of the Chern submanifolds $[m_{\lambda}(\eta)]$ of a stably almost complex manifold $M$ are 
related to the Chern numbers of $M$ through the formula
\begin{equation}
\label{CHGEN}
\langle T_1(\eta)\cdot T_2(\nu),[M]\rangle =S^*_{(y)}(M)+\sum_{\{\lambda||\lambda|>0\}} h_{\lambda}(z)S^*_{(y)}([m_{\lambda}(\eta)]).
\end{equation}
\end{Th}
{\bf Proof.} 
For any stably almost complex vector bundle $\eta$ over a base space $M$ define a power series in $\Lambda[z]\otimes U^*(M)$:
$$
{\cal M}(\eta)=1+\sum_{\{\lambda||\lambda|>0\}} 
h_{\lambda}(z)m_{\lambda}(\eta).
$$
Applying the Riemann-Roch theorem to the constant map $p:M\to \{pt\}$, 
the element ${\cal M}(\eta)\in \Lambda[z]\otimes U^*(M)$ and 
the transformation 
$$
S^*_{(y)}\circ ch_U: U^*(\cdot)\otimes \Lambda[z]\to H^*(\cdot,\Lambda[y]\otimes\Lambda[z]\otimes {\mathbb Q}),
$$
we obtain
\begin{equation}
\label{rmgen}
p^H_!\left(\left(S^*_{(y)}\circ ch_U({\cal M}(\eta))\right)\cdot T(\nu(M))\right)= S^*_{(y)}\circ ch_U\left(p^U_!{\cal M}(\eta)\right),
\end{equation}
where $T(\xi)$ is the Todd genus of the transformation $ S^*_{(y)}\circ ch_U$.
It is obvious that the right hand side of~(\ref{rmgen}) coincides with the 
right hand side of~(\ref{CHGEN}).
It is also a simple corollary of~(\ref{T2}) that $T(\xi)=T_2(\xi)$.
By the splitting principle to calculate $ S^*_{(y)}\circ ch_U({\cal M}(\eta))$ it is sufficient to deal with only 
one-dimensional vector bundles, i.e. with the tautological bundle over
$CP(\infty)$. In this case 
$$
{\cal M}(\eta_1)=\prod^{\infty}_{i=1}(1-z_iu)^{-1},
$$
where $u$ is the canonical generator of $U^*(CP(\infty))$. Applying~(\ref{T2})
we obtain that the right hand side of~(\ref{rmgen}) equals 
$\langle T_1(\eta)\cdot T_2(\nu),[M]\rangle$. 

\

\noindent
Various divisibility conditions, which are  analogues of those~\cite{A} 
coming from  
Atiyah--Singer theorem, can be extracted from Theorem~\ref{ChCh}. 
Consider the genera
$$
T_{A}(\eta)=\prod^m_{i=1} 
\prod^{\infty}_{k=1}\left(1-\frac{z_k(1-{\mathrm e}^{-x_i})}
{\prod^{\infty}_{j=1}\left(1-y_j(1-{\mathrm e}^{-x_i})\right)}\right)^{-1},
$$
$$
T_B(\eta)=
\prod^m_{i=1}\left(\frac{1-{\mathrm e}^{-x_i}}{x_i}
\prod^{\infty}_{j=1}\left(1-y_j(1-{\mathrm e}^{-x_i})\right)^{-1}\right)
$$
where $x_i$, $i=1,...m$, are the 
Chern roots in integral cohomology of a complex vector bundle $\eta$, 
$\dim_{\mathbb C}\eta=m$. 
\begin{Cor}
\label{AS}
For a complex vector bundle $\eta$ over a stably almost complex
manifold $M^{2n}$ with the normal bundle $\nu$, the 
Chern class $T_A(\eta)T_B(\nu)\in H^*(M^{2n},\Lambda[y]\otimes\Lambda[z]\otimes \mathbb Q)$ evaluated on the fundamental 
cycle $[M^{2n}]$
lies in $\Lambda[y]\otimes\Lambda[z]$.
\end{Cor}
{\bf Proof.}  Consider variables 
$\tilde y_j=\tilde y_j(y)$, $j=1,2,...$, such that
$$
S^*_{(\tilde y)}(x)=\prod^{\infty}_{j=1}(1-\tilde y_j x)^{-1}=
\frac{1-{\mathrm e}^{-x}}{x}\prod^{\infty}_{j=1}\left(1-y_j(1-{\mathrm e}^{-x})\right)^{-1}.
$$
Then for any stably almost complex manifold $L$ we have that $S^*_{(\tilde y)}(L)$ is a combination of 
$K$-theory characteristic numbers with coefficients in 
$\Lambda[y]$~\cite{Sm}. 
According to the Hattori--Stong 
theorem~\cite{Ht,St} every $K$-theory characteristic number is integral.
Thus, the corollary follows from Theorem~$\ref{ChCh}$ after making the
substitution 
$y\to\tilde y$.

\begin{Cor}
Let $M^{2n}$ be an almost complex manifold whose tangent space contains 
a complex line subbundle. Then the Euler characteristic of $M^{2n}$ is even.
\end{Cor}
{\bf Proof.} Denote the normal bundle of $M^{2n}$ by $\nu$. Let $\tau\cong 
\tau_1\oplus \eta_1$ be the splitting of the tangent bundle 
$\tau\cong\tau(M^{2n})$, where $\eta_1$ is  a complex line bundle.
From Corollary~\ref{AS} we deduce that the power series 
$\langle T_A(\eta_1\oplus \nu) T_B(\nu),[M^{2n}]\rangle$ has integral 
coefficients. Denote the Chern roots in ordinary cohomology
of $\tau_1$ by $x_1,x_2,...,x_{n-1}$,
and the first Chern class in ordinary cohomology of $\eta_1$ by $x$.
Substitute $z_1=z$, $z_2=z_3=...=0$ and $y_1=y_2=...=0$ into $T_A(\tau_1)$ and
$T_B(\nu)$. Because $T_A(\eta_1\oplus \nu)=T^{-1}_A(\tau_1)$, 
$T_B(\nu)=T^{-1}_B(\tau_1)T^{-1}_B(\eta_1)$  we obtain, 
that the coefficient of $z^{n-1}$ in the polynomial
$$
\left<\left(\frac{x}{1-{\mathrm e}^{-x}}
\prod^{n-1}_{i=1} 
\left(\frac{x_i(1-z(1-{\mathrm e}^{-x_i}))}
{1-{\mathrm e}^{-x_i}}\right)\right) 
,[M^{2n}]\right>
$$
is an integer. One can easily calculate that this coefficient is 
$$
\left<\frac{x}{2}\prod^{n-1}_{i=1}x_i, [M^{2n}]\right>=\frac{\chi(M^{2n})}
{2},
$$
where $\chi(M^{2n})$ is the Euler characteristic of $M^{2n}$.
Thus, $\chi(M^{2n})$ is even. 

\section*{Conclusion}
 
We would like to point out that there is a real version of  
most of the material of the paper. To reproduce the results above in
oriented cobordism theory one notes that the forgetful map 
$\Omega^U_*\to \Omega^{SO}_*/{\mathrm {Tors}}$ is an epimorphism.
There is also the version of the paper 
 for Stiefel--Whitney classes.
It will be interesting to understand what both the versions mean in terms of 
triangulated manifolds.

 The final remark concerning Theorem~\ref{relations} is that  it
gives a complete description of the Landweber--Novikov algebra action 
on the cobordism ring of a point. This might be helpful for calculations of 
homotopy groups of spheres.

\

\noindent
{\bf Acknowledgments.} I express deep gratitude to 
Prof. M. F. Atiyah, Prof. V. M. Buchstaber, Prof. E. G. Rees  
for numerous helpful discussions and advice.
This work was supported by a Seggie Brown Fellowship.

\end{document}